\documentclass[a4paper,10pt]{article}
\usepackage{amsfonts}
\usepackage[latin1]{inputenc}
\usepackage{amssymb , amsmath, amsthm, amsopn, amstext, amscd}
\usepackage[T1]{fontenc}
\usepackage{latexsym, enumerate}
\usepackage{graphicx}

\newtheorem{prp}{Proposition}
\newtheorem{defi}{Definition}
\newtheorem{thm}{Theorem}
\newtheorem{lm}{Lemma}
\newtheorem{rmq}{Remark}

\title{Nonparametric estimation of the purity of a quantum state in quantum homodyne tomography with noisy data.}

\author{Katia M\'{e}ziani}

\begin{document}

\maketitle

\begin{abstract}
The aim of this work  is to estimate a quadratic functional of a unknown Wigner function from noisy tomographic data. The Wigner function can be seen as the representation of the quantum state of a light beam. The estimation of a quadratic functional is done from result of quantum homodyne measurement performed on identically prepared quantum systems.\\
We start by constructing an estimator of a quadratic functional of the Wigner function. We show that the proposed estimator is optimal or nearly optimal in a minimax sense over a class of infinitely differentiable functions. Parametric rates are also reached for some values of the smoothness parameters and the asymptotic normality is given. Then, we construct an adaptive estimator that does not depend on the smoothness parameters and prove it is minimax over some set-ups.
\end{abstract}
\textbf{AMS 2000 subject classifications:} 62G05, 62G20, 81V80,\\
\textbf{Key Words:} Adaptive estimation, deconvolution, infinitely differentiable functions, minimax risk, quadratic functional estimation, quantum state, Wigner function, Radon transform, quantum homodyne tomography, asymptotic normality.

\section{Introduction}
\label{sec:1}
In quantum mechanics, the quantum state of a system completely describes all aspects of the system. The instantaneous state of a quantum system encodes the probabilities of its measurable properties, or "observables" (examples of observables include energy, position, momentum and angular momentum). Generally, quantum mechanics do not assign determinist values to 
observables. Instead, it makes predictions about probability distributions; that is, the probability of obtaining each of the possible outcomes from measuring an observable.\\
We have two mathematical representations of a quantum state: the density matrix $\rho$ and its associated Wigner function $W_\rho$. The densitymatrix $\rho$, which describes completely a quantum state, is hermitian, positive definite and with trace one. It can be finite or infinite dimensional. Equivalently the corresponding Wigner function $W_\rho:\mathbb{R}^2\rightarrow\mathbb{R}$ may be defined. In general, $W_\rho$ is regarded as a generalized probability density, integrating to plus one over the whole plane. It does not satisfy all the properties of a proper probability density as it can, and normally does, go negative for states which have no classical model. It satisfies also certain intrinsic positivity constraints in the sense that it corresponds to a density matrix.
\paragraph{}
In this paper we address the  problem of estimating the quadratic functional $d^2=\int W_\rho^2$ of the Wigner function of a monochromatic light in a cavity prepared in the state $\rho$ by using Quantum Homodyne Tomographic (QHT\footnote{We refer the interested reader to Artiles \textit{et al.} (2005) \cite{MR2136642} for further details on the physical background}) 
data measurement performed on independent, identical systems. The Quantum Homodyne Detection (QHD) has been put in practice for the first time by Smithey \textit{et al.} (1993) \cite{SmitA}, we will detail this technique in section~\ref{sec:2.2}.\\
We study the quantity $d^2=\int W_\rho^2$ which has an interest in itself as a physical measure of the purity of quantum state. It allows us to detect pure state and mixed state as it always equals $\frac{1}{2\pi}$ in case of pure states. A state is called pure if it cannot be represented as a mixture (convex combination) of other states, i.e., if it is an extreme point of the convex set of states. All other states are called mixed states.\\
The QHD technique gives results of the measure of the electric and the magnetic fields $(p,q)$ of the studied laser for some phase $\Phi$. In the ideal case, we would observe the random variable $(X,\Phi)=(\cos(\Phi)Q+\sin(\Phi)P,\Phi)$ where $\Phi$ is chosen independently of $(Q,P)$, and uniformly in the interval $[0,\pi]$. In our paper we do not consider the ideal data $ (X,\Phi)$ but the noisy observations $(Y,\Phi)$ where $Y$ is the sum of the random variable $X$ and a gaussian random variable $\xi$. We assume that the unknown function $W_\rho$ belongs to $\mathcal{A}(\alpha,r,L)$ a class of super smooth functions where $\alpha>0$, $0<r\leq 2$ and $L>0$ will be defined later. Those classes are similar to those of Cavalier (2000) \cite{MR1790012} for $r=1$ and functions are defined on $\mathbb{R}^d$; Butucea and Tsybakov (2007) \cite{ButuTsyba04} on $\mathbb{R}$; Butucea \textit{et al.} (2007) \cite{ButGutArt05} on $\mathbb{R}^2$.\paragraph{}
The study of quadratic functionals started with Bickel and Ritov (1988) \cite{BickelRitov}, who have considered the problem of estimating the integral of the square of a derivative of a probability density function and obtained nonparametric rates. Their results have been extended by Birg\'{e} and Massart (1995) \cite{MR1331653} on the estimation of more general functionals, who established nonparametric lower bounds. The study of general functionals was completed by Kerkyacharian and Picard (1996) \cite{MR1394973} for minimax rates. Laurent (1996) \cite{MR1394981} gave efficient estimation of some functionals of a density function at parametric rate. The problem of adaptive estimation of general functionals has been 
considered by Tribouley (2000) \cite{MR1772223} in the classical white noise model.\\
In the convolution model, Butucea (2004) \cite{Butu04} has estimated a quadratic functional of a density on $\mathbb{R}$ and applied it to the goodness-of-fit test in $L_2$ distance.\\
In our paper, the first difficulty is that we do not deal with proper probability density function but with quasiprobability density. Moreover, note that our problem is a double inverse problem as we observe the Radon transform of $W_\rho$ (PET) with a convolution (white noise).\\
Inverse problems have been extensively studied in mathematical literature. In a positron emission tomography (PET) perspective, the problem of estimating a probability density on $\mathbb{R}^2$ from tomographic data $(X_k,\Phi_k)$ has been treated by Korostel\"{e}v and Tsybakov (1993) \cite{MR1226450} and johnstone and Silverman (1990) \cite{MR1041393}. Cavalier (2000) \cite{MR1790012} considered also PET model and obtained an estimator of a multi-dimensional density function which is asymptotically sharp minimax, i.e. it achieves the optimal rate of convergence and attains the best constant for the minimax risk.\\
The estimation of the Wigner function $W_\rho$ has been treated by Gu{\c{t}}{\u{a}} and Artiles (2006) \cite{GutaArt05} in the case free of noise. Our noisy model has been studied in a parametric framework by D'Ariano and in a nonparametric framework for the estimation of the Wigner function by Butucea \textit{et al.} (2007) \cite{ButGutArt05}. We propose to estimate the integral of the square  of the Wigner function rather than the function itself.\\
Other problems have been considered, in the context of tomography: Goldenshluger and spokoiny (2006) \cite{GoldSpok06} have considered the problem of recovering edges of an image from noisy tomographic data in a white noise model and reached nearly optimal rate. Recovering boundaries in models that involve indirect observations in the $d$-dimensional Euclidean space $\mathbb{R}^d$ has been discussed recently in Goldenshluger and Zeevi (2006) \cite{GoldZeev06}. We note that a Wigner function cannot have a bounded support.
\paragraph{}
The main contributions of this paper are the following. We propose a method for estimating a quadratic functional of a generalized probability density which may take negative values from indirect and noisy observations in view to detect pure states and mixed states. It is shown that the proposed estimator is optimal or nearly optimal in a minimax sense -depending on the smoothness parameter $r$ of the class $\mathcal{A}(\alpha,r,L)$. Moreover, an adaptive estimator is constructed which attains optimal rates. Another main interest of the estimation of $d^2$ is the important application to goodness-of-fit test in $\mathbb{L}_2$-norm in quantum statistics. This means that physicists want to test whether they produced a laser in the quantum state $\rho_0$ or something different. This can be done via the Wigner functions as follows:
\[\left\lbrace
  \begin{array}{c l}
    H_0: & \text{ $W_\rho=W_{\rho_0}$},\\
    H_1: &
\text{$\sup_{W_\rho\in\mathcal{A}(\alpha,r,L)}\|W_\rho-W_{\rho_0}\|_2\geq
c\cdot\varphi_n$}.
  \end{array}
\right. \] where $\varphi_n$ is a sequence which tends to 0 when $n\rightarrow\infty$ and it is the testing rate. We can device a test statistic based on the estimator of $d^2=\int W_\rho^2$ constructed in this paper. Similary to Butucea (2004) \cite{Butu04} we conjecture that the testing rates are  of the same order as the nonparametric ones found in this paper.
\paragraph{}
The rest of the paper is organized as follows. In Section~\ref{sec:3} we formulate the statistical model and introduce notation and properties of quantities of interest. In Section~\ref{sec:4} we construct an estimator of the quadratic functional of the unknown Wigner function, along with the bias-variance decomposition. Our main theoretical results are presented in Section~\ref{sec:5}. In Section~\ref{sec:6}, we derive some example of quantum states. Proofs of upper and lower bounds are given in Sections~\ref{sec:7},~\ref{sec:8} and ~\ref{sec:9}.

\section{Preliminaries}
\label{sec:2}

\subsection{Definition}
\label{sec:2.1}

We study this problem in a minimax framework. Let $d_n^{2}$ be an estimator of $d^{2}=\int{W_\rho^2}$ based on this indirect noisy observations $(Y_i,\Phi_i)$, $i=1,\ldots,n$ as anounced above. We measure the accuracy of $d_n^{2}$ by the maximal risk$$\mathcal{R}(d_n^{2};\mathcal{A}(\alpha,r,L))=\sup_{W_\rho\in\mathcal{A}(\alpha,r,L)}E_{\rho}[|d_n^2-d^2|^2]$$
over the class $\mathcal{A}(\alpha,r,L)$. Here $E_{\rho}$, $P_{\rho}$ denote the expected value and probability when the true underlying quantum state is $\rho$. The minimax risk is defined by $$\mathcal{R}^{*}(\mathcal{A}(\alpha,r,L))=\inf_{\widehat{d}_n^2}\mathcal{R}(\widehat{d_n^{2}};\mathcal{A}(\alpha,r,L))$$
where the infimum is taken over all possible estimators $\widehat{d}_n^2$ of the quadratic functional of the Wigner function $W_{\rho}$.\\
Let $\varphi_n$ be a positive sequence, an estimator $d^2_n$ is \textit{optimal in a minimax sense}
\begin{itemize}
\item if it satisfies the following \textit{upper bound}
\begin{eqnarray}
\label{UB}
\limsup_{n\to\infty}\varphi_n^{-2}\mathcal{R}(d_n^{2};\mathcal{A}(\alpha,r,L))\leq C_{u},
\end{eqnarray}
\item and if the following \textit{lower bound} is satisfied
\begin{eqnarray}
\label{LB}
\liminf_{n\rightarrow\infty}\inf_{\widehat{d}_n^2}\varphi_n^{-2}\mathcal{R}(\widehat{d}_n^2;\mathcal{A}(\alpha,r,L))\geq C_{l},
\end{eqnarray}
\end{itemize}
where the infimum is taken over all possible estimators $\widehat{d}_n^2$ of the quadratic functional of the Wigner function $W_{\rho}$.Then, $\varphi_n$ is called \textit{optimal rate in a minimax sense}. Our aim is to find rate optimal estimator of $d^{2}$ and to establish asymptotics of minimax risks for some classes of Wigner functions $\mathcal{A}(\alpha,r,L)$. We  rely on Butucea \textit{et al.} (2007) \cite{ButGutArt05}, who derived rate optimal pointwise and adaptive estimators of $W_\rho$ (instead of $\int{W_\rho^2}$ in our case) from indirect noisy observations.

\subsection{Quantum Homodyne Tomography}
\label{sec:2.2}

The theoretical foundation of quantum state reconstruction was outlined by Vogel and Risken (1989) \cite{Vogelrisken} and has inspired the first experiments determining the quantum state of a light field, initially with optical pulses with Smithey \textit{et al.} (1993) \cite{SmitA} and Smithey \textit{et al.} \cite{Smitheybis}.\\
\begin{figure}[htbp!]
\begin{center}
\includegraphics[trim = 0cm 3cm 0mm 0cm,clip,width=8cm]{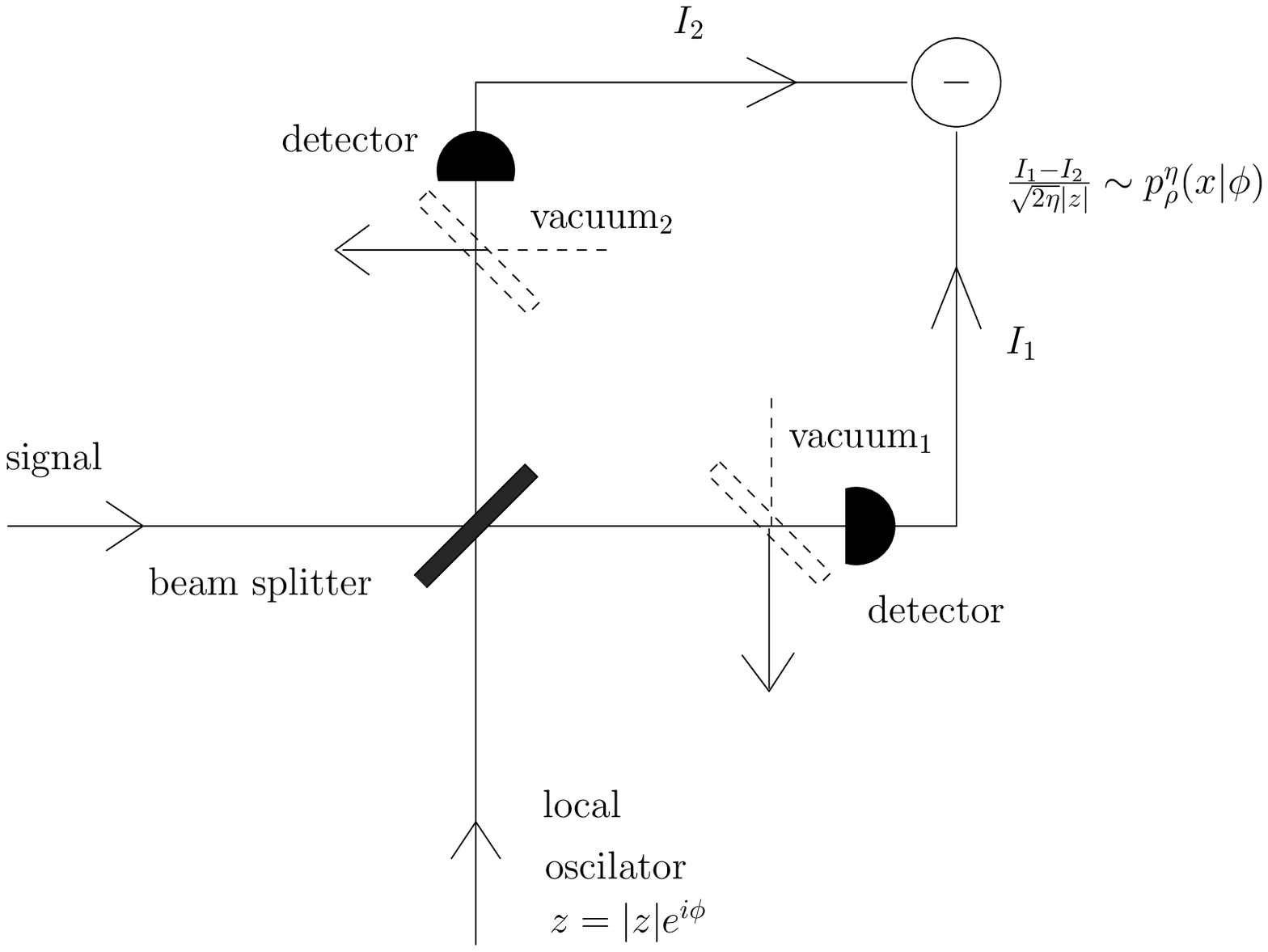}
\caption{QHT mesurement}
\label{fig:1}
\end{center}
\end{figure}
\\
The physicists developed a monochromatic laser in state $\rho$ in a cavity. In order to study it, one takes measurement by quantum tomography homodyne (QHT). This technique schematized in figure~\ref{fig:1} consists in mixing the laser to be studied with a laser of reference of high intensity $\left|z\right|>>1$ called local oscillator. Then the beam obtained is split into two and two photodetectors measure each one of the beams ($I_1,I_2$). One measures $X$ the difference of the intensities of the two beams and rescale it by the intensity $\left|z\right|$. Thus for the cavity pulse chosen to be $\phi$, data $(X,\Phi)$ should be obtained. It is widely known in the physical litterature (see Leonhardt (1997) \cite{Leon97}) that an additive gaussian noise is mixed with ideal data $X$, giving for known efficiency $\eta$, data $Y$.

\section{Statistical context}
\label{sec:3}

\subsection{Problem formulation}
\label{sec:3.1}

In the present paper we estimate the integral of the square of the Wigner function from data measurement performed on $n$ identical quantum systems where the Wigner function is assumed to be a joint generalized density of two variables $P$ and $Q$, $W_\rho:\mathbb{R}^2\rightarrow\mathbb{R}.$ It may take negative values but it integrates to plus one over the whole plane. For further information of the Wigner function, we invite readers to refer to the paper by Artiles \textit{et al.} \cite{MR2136642}.\\
Our statistical problem can been seen as follow: consider $(X_1,\Phi_1)\ldots(X_n,\Phi_n)$ independent identically distributed random 
variables with values in $\mathbb{R}\times[0,\pi]$. The probability density of $(X,\Phi)$ equals the \textbf{Radon transform} \textbf{$\Re[W_\rho]$} of the Wigner function with respect to the measure $\lambda/\pi$, where $\lambda$ is the Lebesgue measure on $\mathbb{R}\times[0,\pi]$.

\begin{equation}
\label{radon} p_\rho(x/\phi):=\Re[W_\rho](x,\phi)=\int_{-\infty}^\infty W_\rho(x\cos\phi+t\sin\phi,x\sin\phi-t\cos\phi)dt
\end{equation}and $X$ has density $p_\rho(x/\phi)$. As we annouced in the introduction we do not observe the ideal data $(X_\ell,\Phi_\ell)$ $\ell=1,\ldots,n$ but a degraded noisy version $(Y_1,\Phi_1)\ldots(Y_n,\Phi_n)$,
\begin{equation}
\label{model} Y_\ell:=\sqrt{\eta}X_\ell+\sqrt{(1-\eta)/2}\xi_\ell
\end{equation}with $\xi_\ell$ a standard Gaussian random variables independent of all $(X_k,\Phi_k)$ and $0<\eta<1$ is a known parameter. The parameter $\eta$ is called the detection efficiency and represents the proportion of photons which are not detected due to various losses in the measurement process. We note $p_\rho^\eta(x,\phi)$ the density of $(Y_\ell,\Phi_\ell)$. Thus, $p_\rho^\eta(x,\phi)$ is the convolution of the density $\frac{1}{\sqrt{\eta}}p_\rho^\eta(\frac{x}{\sqrt{\eta}},\phi)$ with the density of a centered Gaussian density having variance $(1-\eta)/2$. We assume that the unknown Wigner function $W_\rho$ belongs to a class $\mathcal{A}(\alpha,r,L)$ of infinitely differentiable functions. For $0<r\leq 2$, $\alpha>0$ and $L>0$ define
\begin{equation}
\label{ens fctnel} \mathcal{A}(\alpha,r,L)=\{W_\rho:\int_{\mathbb{R}^{^2}}|\widetilde{W_\rho}(u,v)|^{2}e^{2\alpha\|(u,v)\|_2^r} du
dv\leqslant(2\pi)^{2}L\}
\end{equation}where $\|(u,v)\|_2=\sqrt{u^2+v^2}$ is the euclidian norm.

\subsection{Properties of Wigner functions and remarkable equations}
\label{sec:3.2}

In this paragraph we will state some very useful properties the Wigner function.

\paragraph{\textbf{Fourier transforms}}
A remarkable relation links the Fourier transform of the Wigner function to the Fourier transform of its Radon transform. If we denote
\begin{eqnarray*}
\widetilde{W}_{\rho}(u,v)&:=&\mathcal{F}_2[W_\rho](u,v),
\end{eqnarray*}then
\begin{eqnarray}
\label{fourierp}
\widetilde{W}_{\rho}(t\cos\phi,t\sin\phi)&:=&\mathcal{F}_1[p_\rho(\cdot/\phi)](t)=E_\rho[e^{itX}]
\end{eqnarray}where $\mathcal{F}_2$, $\mathcal{F}_1$ denote the fourier transform w.r.t two, respectively one variables.

\paragraph{\textbf{Some remarkable equations}}
In Section~\ref{sec:8}, most of the proofs make extensive use of the following equations. Since

\begin{eqnarray*}
E_\rho[e^{itY}]&=&E_\rho[e^{it\sqrt{\eta}X}]\cdot E_\rho[e^{it\sqrt{\frac{1-\eta}{2}}\xi}]
\end{eqnarray*}then
\begin{eqnarray}
\label{fourierproun}
\mathcal{F}_1[p_\rho^\eta(\cdot/\phi)](t)&=&\mathcal{F}_1[\frac{1}{\eta}p_\rho^\eta(\frac{.}{\eta}/\phi)](t)\cdot\widetilde{N}^\eta(t)\\
\label{fourierprodeux}
&=&\mathcal{F}_1[p_\rho(\cdot/\phi)](\sqrt{\eta}t)\cdot\widetilde{N}^\eta(t),
\end{eqnarray}where $\widetilde{N}^\eta(t)$ denotes the Fourier transform of $\sqrt{(1-\eta)/2}\xi\sim\mathcal{N}(0;(1-\eta)/2).$ Then
\begin{equation}
\label{fourierbruit}
\widetilde{N}^\eta(t):=E_\rho[e^{it\sqrt{(1-\eta)/2}\xi}]=e^{-\frac{t^2}{4}(1-\eta)}.
\end{equation}
\section{Estimation procedure}
\label{sec:4}

We are now able to define the estimation procedure of the quadratic functional $d^2=\int W_\rho^2$ of the unknown function $W_\rho$ directly from data $(Y_{\ell},\phi_{\ell})$. Next we evaluate an upper bound of the maximal risk uniformly over all Wigner functions in the class $\mathcal{A}(\alpha,r,L)$.

\subsection{ Kernel estimator}
\label{sec:4.1}

Let us define our estimator as a U-statistic of order 2:

\begin{defi}
Let $(Y_{\ell},\phi_{\ell}),\ell=1,\ldots,n$, be i.i.d data coming from the model \eqref{model}, and $\delta=\delta_n\rightarrow 0$ as $n\rightarrow\infty$. The estimator $d_n^{2}$ can be written
\begin{equation}
\label{estimateur} 
d_{n}^{2}:=\frac{1}{(2\pi)^{2}}\frac{1}{n(n-1)}\sum_{k\neq\ell=1}^n\int_{|t|\leq\frac{1}{\delta\sqrt{\eta}}}\int_0^\pi\eta|t|e^{\frac{t^2}{2}(1-\eta)}e^{itY_k-itY_\ell}d\phi dt.
\end{equation}
\end{defi}

\begin{defi}
\label{def}
Let $d_n^{2}$ be the estimator defined in \eqref{estimateur}, having bandwidth $\delta>0$. We call the bias and the variance of the estimator, respectively: $$B(d_n^{2}):=|E_\rho[d_n^{2}]-d^{2}|^2\quad \textrm{and}\quad\text{Var}(d_n^{2}):=E_\rho\left[|d_n^{2}-d^{2}|^2\right].$$
\end{defi}

\subsection{Bias-variance decomposition}
\label{sec:4.2}

The following proposition plays an important role in the proof of the upper bound of the risk as we split it into the bias term and the variance term.
\begin{prp}
\label{prop}
Let $(Y_{\ell},\phi_{\ell}),\ell=1,\ldots,n$ be i.i.d data coming from the model \eqref{model} and $d_{n}^{2}$ be the estimator in \eqref{estimateur} (with $\delta\rightarrow 0$ as $n\rightarrow \infty$) of $d^2$ the quadratic functionnal of the Wigner function $W_\rho$ which is lying in the class $\mathcal{A}(\alpha,r,L)$ with $\alpha>0$, $L>0$ and $0<r\leq 2$ defined in \eqref{ens fctnel} then,
\begin{enumerate}
\item for all $0<r\leq 2$
\begin{eqnarray}
\label{propa} |E_\rho[d_n^{2}]-d^{2}|^2 &\leq &L^2 e^{-4\alpha/\delta^r},
\end{eqnarray}
\item for all $0<r<2$
\begin{eqnarray}
\label{propbun} 
\text{Var}(d_n^{2})\leq\frac{8\eta^2/(1-\eta)^2}{\pi^2 n^2}e^{\frac{1-\eta}{\eta}\frac{1}{\delta^2}}+\frac{8L}{n\pi}\frac{\eta}{1-\eta}e^{\frac{1-\eta}{2\eta}\frac{1}{\delta^2}-\frac{2\alpha}{\delta^r}},
\end{eqnarray}
\item for all $r=2$ and $\frac{1-\eta}{2\eta}-2\alpha>0$
\begin{eqnarray}
\label{propbdeux} 
\text{Var}(d_n^{2})\leq\frac{8\eta^2/(1-\eta)^2}{\pi^2
n^2}e^{\frac{1-\eta}{\eta}\frac{1}{\delta^2}}+\frac{8L}{n\pi}\frac{\eta}{1-\eta-4\alpha\eta}e^{(\frac{1-\eta}{2\eta}-2\alpha)\frac{1}{\delta^2}},
\end{eqnarray}
\item for all $r=2$ and $\frac{1-\eta}{2\eta}-2\alpha<0$
\begin{eqnarray}
\label{propbtrois} 
\text{Var}(d_n^{2})\leq\frac{8\eta^2/(1-\eta)^2}{\pi^2
n^2}e^{\frac{1-\eta}{\eta}\frac{1}{\delta^2}}+\frac{1}{n}\cdot\frac{8\eta L}{4\alpha\eta-1+\eta}.
\end{eqnarray}
\end{enumerate}
\end{prp}
The proof of this proposition is given in section~\ref{sec:8}.

\section{Main results}
\label{sec:5}

In this section, the first theorem considers the case of nonparametric rates of convergence of our estimator which is proven optimal or nearly optimal (as we loose a logarithmic factor in the lower bound) in the minimax sense. In the second theorem, our estimator attains the parametric rate $1/n$.
\begin{thm}
\label{theo:1}
Let $(Y_{\ell},\phi_{\ell}),\ell=1,\ldots,n$ be i.i.d data coming from the model \eqref{model} where the underlying parameter is the Wigner function $W_\rho$ lying in the class $\mathcal{A}(\alpha,r,L)$, $\alpha>0$ and $L>0$. Then for $d_{n}^{2}$ defined in \eqref{estimateur} and according to the definition given to section~\ref{sec:2.1},
\begin{enumerate}
\item for $0<r<2$, with $\delta:=\delta_{opt}$ solution of the equation
\begin{eqnarray}
\label{ftre}
\frac{1-\eta}{2\eta}\frac{1}{\delta_{opt}^2}+\frac{2\alpha}{\delta_{opt}^r}=\log
n-(\log\log n)^2,
\end{eqnarray}
we reach the optimal rate $\varphi_n$ with $C_u=1$, $C_l=1/16$ constants defined in \eqref{UB} and \eqref{LB}
\begin{eqnarray}
\label{vitesse}
\varphi_n^2=L^2e^{\frac{-4\alpha}{\delta_{opt}^r}},
\end{eqnarray}
\item for $r=2$, $\frac{1-\eta}{2\eta}-2\alpha>0$ and by taking $\delta=\delta^{*}=\left(\frac{\log 
n}{\frac{1-\eta}{2\eta}+2\alpha}\right)^{-1/2}$, the rate of convergence is nearly optimal as
\begin{eqnarray}
\label{vitesse1}
\varphi_n^2=n^{\frac{-4\alpha}{\frac{1-\eta}{2\eta}+2\alpha}},
\end{eqnarray}  is the rate of convergence in the upper bound \eqref{UB} and
\begin{eqnarray}
\label{vitesse2}
\varphi_n^2=(n\log n )^{\frac{-4\alpha}{\frac{1-\eta}{2\eta}+2\alpha}},
\end{eqnarray} is the rate of convergence in the lower bound \eqref{LB}.
\end{enumerate}
\end{thm}
To prove the Theorem~\ref{theo:1}, one has to prove on the one hand, the upper bound (section~\ref{sec:7}) and on the other hand, the lower bound (section~\ref{sec:9}) according to the definition given in section~\ref{sec:2.1}.
\begin{thm}
\label{theo:2}
Let $(Y_{\ell},\phi_{\ell}),\ell=1,\ldots,n$ be i.i.d data coming from the model \eqref{model} where the underlying parameter is the Wigner function $W_\rho$ lying in the class $\mathcal{A}(\alpha,r,L)$, $r=2$, $\alpha>0$, $L>0$ and $\frac{1-\eta}{2\eta}-2\alpha<0$. Then for $d_{n}^{2}$ defined in \eqref{estimateur} with $\delta=\delta^*=\left(\frac{\eta\log n}{1-\eta}\right)^{-1/2}$, the rate of convergence is parametric: $\varphi_n^2=\frac{1}{n}$.\\
Moreover, in this case our estimator \eqref{ens fctnel} is asymptotically normally distributed $$\sqrt{n}(d^2_n-d^2)\rightarrow\mathcal{N}(0,\mathcal{W}),$$ with asymptotic variance
$$\mathcal{W}=\frac{1}{4\pi^{2}}\int\int|t_1||t_2|e^{\frac{1-\eta}{2\eta}t_1t_2}E[e^{it_1X}]E[e^{it_2X}] E[e^{-i(t_1+t_2)X}] dt_1 dt_2-4d^2.$$
\end{thm}
The proof of the Theorem~\ref{theo:2} in given in section~\ref{sec:7}.
\begin{rmq}
\label{rmq:1}
We are able to give a more explicit form for the bandwidth and thus for the bias term which is asymptotically equivalent to the rate according to the values of $r$. Let $s_n:=\frac{\log n-(\log\log n)^2}{2a}$ where $a:=\frac{1-\eta}{4\eta}$, then we make successive approximations in \eqref{ftre} starting with $\delta_0$ and we plug it back into \eqref{ftre}, we find $\delta_1$. And successively for all $k\geq 1$ we have $\delta_k$. Values are given by table~\ref{tab:1} and table~\ref{tab:2}.\\
\begin{table}[htbp!]
\caption{Procedure}
\label{tab:1}
\begin{tabular}{lll}
\hline\noalign{\smallskip}
$\delta_0$ & $\delta_1$ & for all $k\geq 1$, $\delta_k$ \\
\noalign{\smallskip}\hline\noalign{\smallskip}
$=s_n^{-1/2}$ & $=(s_n-\frac{\alpha}{a}\delta_0^{-r})^{-1/2}$ & 
$=(s_n-\frac{\alpha}{a}\delta_{k-1}^{-r})^{-1/2}$ \\
\noalign{\smallskip}\hline
\end{tabular}
\end{table}
\begin{table}[htbp!]
\caption{Rates of convergence}
\label{tab:2}
\begin{tabular}{lll}
\hline\noalign{\smallskip}
If $r$  & It is enough  & And the   \\
belongs to & to choose & rate is  \\
\noalign{\smallskip}\hline\noalign{\smallskip}
$r\in]0,1]$ & $\delta=\delta_1$ & $L^2e^{\left(-4\alpha s_n^{r/2}+o(1)\right)}$ \\
$r\in]1,4/3]$ & $\delta=\delta_2$ & $L^2e^{\left(-4\alpha s_n^{r/2}+C_1 s_n^{r-1}-o(1)\right)}$ \\
$r\in]\frac{2(k-1)}{k},\frac{2k}{k+1}]$ & $\delta=\delta_k$ & $L^2 e^{\left(-4\alpha s_n^{r/2}+C_1
s_n^{r-1}-\ldots+C_{k-1}s_n^{kr/2-(k-1)}+o(1)\right)}$ \\
\noalign{\smallskip}\hline
\end{tabular}
\end{table}
\end{rmq}
In the previous theorems, the bandwidth $\delta_{opt}$ depends on the parameters $\alpha$, and $r$ of the class $\mathcal{A}(\alpha,r,L)$ which may be difficult to evalute in practice. However, it is possible to construct an adaptive estimator which does not depend on these parameters and which attains the same asymptotic behavior as in Theorem~\ref{theo:1}, provided that these parameters lie in a certain set. Note that the parameter $\eta$ is supposed to be known. Define two sets of parameters
\begin{eqnarray*}
\Theta_1&=&\{(\alpha,r,L):\alpha>0,L>0,0<r<1\}\\
\Theta_2&=&\{(\alpha,r,L):0<\alpha\leq\alpha_0,L>0,r=1\},\quad\alpha_0>0.
\end{eqnarray*}
\begin{thm}
\label{theo:3}
Let $(Y_{\ell},\phi_{\ell}),\ell=1,\ldots,n$ be i.i.d data coming from the model \eqref{model}. For $\delta=\delta^i_{ad},\,i=1,2$, let $d_{\delta,n}^{2}$  be the estimator defined by
$$d_{\delta,n}^{2}:=\frac{1}{(2\pi)^{2}}\frac{1}{n(n-1)}\sum_{k\neq\ell=1}^n\int_{|t|\leq\frac{1}{\delta^i_{ad}\sqrt{\eta}}}\int_0^\pi\eta|t|e^{\frac{t^2}{2}(1-\eta)}e^{itY_k-itY_\ell}d\phi dt,$$
with $\delta^1_{ad}=(\frac{2 \eta \log n}{1-\eta}-\sqrt{\frac{2\eta\log n}{1-\eta}})^{-1/2}$ and $\delta^2_{ad}=(\frac{2\eta\log n}{1-\eta}-\frac{4A\eta}{1-\eta}\sqrt{\frac{2\eta\log n}{1-\eta}})^{-1/2}$, $A>\alpha_0$. Then, for all  $(\alpha,r,L)\in\Theta_i,$ $i=1,2$, respectively,
$$\limsup_{n\rightarrow\infty}\sup_{W_\rho\in\mathcal{A}(\alpha,r,L)}E[|d_{\delta,n}^2-d^2|^2]\varphi_n^{-2}\leq C_i,$$ where $\varphi_n^{-2}$ is the rate defined in \eqref{vitesse} and the constants are respectively $C_1=1$ and $C_2=\exp{(\frac{8A\alpha\eta}{1-\eta}-\frac{8\alpha^2\eta}{1-\eta})}.$
\end{thm}
The proof of the adaptive case in given section~\ref{sec:7}.

\section{Examples}
\label{sec:6}

The Table~\ref{tab:3} shows five examples of quantum pure states and one example of mixed state which can be created at this moment in laboratory. Among the pure states we consider the vacuum state which is the pure state of zero photons, the single photon state, the coherent state which characterizes the laser pulse with an average of $N$ photons. The squeezed states have Gaussian Wigner functions whose variances in the two directions have a fixed product. And the well-known Schr\"{o}dinger's Cat which is also a pure state.\\
Note that for pure states, $d^2=1/2\pi$. The thermal state is a mixed state describing equilibrium at temperature equal to $1/\beta$, having Gaussian Wigner function with variance increasing with the temperature. This state is mixed and here we find $d^2=\frac{\tanh(\beta/2)}{2\pi}$. For these examples of quantum states, the procedure gives fast parametric rates with $r=2$ and $\frac{1-\eta}{2\eta}-2\alpha<0$. We can easily check that each Wigner function belongs to the class $\mathcal{A}(\alpha,2,L)$ for small enough values of $\alpha$ (see Table~\ref{tab:3}).
\begin{table}[htbp!]
\caption{Examples of quantum states}
\label{tab:3}
\begin{tabular}{lll}
\hline\noalign{\smallskip}
State & Fourier transform of Wigner & in the class  \\
      & function $\widetilde{W_\rho}(u,v)$ & $\mathcal{A}(\alpha,2,L)$ if  \\
\noalign{\smallskip}\hline\noalign{\smallskip}
Vacuum state & $\exp\left(\frac{-\|(u,v)\|^2_2}{4}\right)$ & $\alpha<1/4$ \\
Single photon state & $\left(1-\frac{\|(u,v)\|^2_2}{2}\right)\exp\left(\frac{-\|(u,v)\|^2_2}{4}\right)$ & $\alpha<1/4$ \\
Schr\"{o}dinger's Cat $X_0>0$ & $\frac{e^{\frac{-\|(u,v)\|^2_2}{4}}}{2(1+e^{-X_0^2})}\left(\cos(2uX_0)+e^{-X_0^2}\cosh (X_0v)\right)$ & $\alpha<1/4$ \\
Coherent state $N\in\mathbb{R}_+$ & $\exp\left(\frac{-\|(u,v)\|^2_2}{4}+i\sqrt{N}v\right)$ & $\alpha<1/4$ \\
Squeezed state $N\in\mathbb{R}_+$, $\xi\in\mathbb{R}$ & $\exp\left(-\frac{u^2}{4}e^{2\xi}-\frac{v^2}{4}e^{-2\xi}+iv\alpha\right)$ & $\alpha<e^{2\xi}/4$ \\
Thermal state $\beta>0$ & $\exp\left(\frac{-\|(u,v)\|^2_2}{4(\tanh(\beta/2))^2}\right)$ & $\alpha<\frac{(\tanh(\beta/2))^2}{4}$ \\ \noalign{\smallskip}\hline
\end{tabular}
\end{table}
\\
Our previous results show that our estimator of the purity atteins the parametric rate $1/n$ if $\eta>\frac{1}{1+4\alpha}$. This is not restrictive at all. In practice, physicists usually find $\eta>0.8$ and more often $\eta$ is close to $0.9$ and $0.95$. Thus, by choosing $\alpha$ as close to its upper bound (in Table~\ref{tab:3}) as possible we make sure that our 
estimator attains the parametric rate.

\section{Proof of the upper bounds of theorems}
\label{sec:7}

\paragraph{\textbf{Sketch of proof of upper bound in Theorem~\ref{theo:1}-\eqref{vitesse}}}
For $0<r<2$ and by \eqref{propa} and \eqref{propbun}
\begin{eqnarray*}
\text{Var}(d_n^{2})&\leq&\frac{8\eta^2/(1-\eta)^2}{\pi^2 n^2}e^{\frac{1-\eta}{\eta}\frac{1}{\delta^2}}+\frac{8L}{n\pi}\frac{\eta}{1-\eta}e^{\frac{1-\eta}{2\eta}\frac{1}{\delta^2}-\frac{2\alpha}{\delta^r}}\\
&=&\frac{C_{V1}}{n^2}e^{\frac{1-\eta}{\eta}\frac{1}{\delta^2}}+\frac{C_{V2}}{n}e^{\frac{1-\eta}{2\eta}\frac{1}{\delta^2}-\frac{2\alpha}{\delta^r}}.
\end{eqnarray*}
On the one hand, we select the bandwidth $\delta^*$ as
$$\delta^*=\arg\inf_{\delta>0}\left\{\frac{C_{V1}}{n^2}e^{\frac{1-\eta}{\eta}\frac{1}{\delta^2}}+\frac{C_{V2}}{n}e^{\frac{1-\eta}{2\eta}\frac{1}{\delta^2}-\frac{2\alpha}{\delta^r}}+C_Be^{-4\alpha/\delta^r}\right\},$$ by taking derivatives, $\delta^*$ is a positive real number satisfying $$\frac{1-\eta}{2\eta}\frac{1}{\delta^{*2}}+\frac{2\alpha}{\delta^{*r}}+\log(\delta^{*r-2})=\log n$$ and we notice that $B(d^2_n)\sim\delta^{r-2}Var(d^2_n)$. So the rate of convergence for the upper bound is given by the bias i.e. $\varphi^2_n=B(d^2_n)(1+o(1))$. On the other hand, we show that by taking $\delta:=\delta_{opt}$ the unique solution of the equation $$\frac{1-\eta}{2\eta}\frac{1}{\delta_{opt}^2}+\frac{2\alpha}{\delta_{opt}^r}=\log n-(\log\log n)^2$$ we obtain the same results. We find $B(d^2_n)\sim\delta^{r-2}Var(d^2_n)$ for $\delta^*$.
\begin{eqnarray*}
\frac{\delta_{opt}^{r-2}}{n}\exp{\left(\frac{1-\eta}{2\eta\delta_{opt}^2}\frac{-2\alpha}{\delta_{opt}^r}\right)}
&=& \frac{\delta_{opt}^{r-2}}{n}\exp{\left(\log n-(\log\log n)^2-\frac{4\alpha}{\delta_{opt}^r}\right)}\\
&=&\frac{\delta_{opt}^{r-2}}{(\log\log n)^2}\exp{\left(\frac{-4\alpha}{\delta_{opt}^r}\right)}\\
&=&\frac{\left(\log n/(2\beta)\right)^{(2-r)/2}}{(\log\log n)^2}\exp{\left(\frac{-4\alpha}{\delta_{opt}^r}\right)}\\
&=&o(1)\exp{\left(\frac{-4\alpha}{\delta_{opt}^r}\right)}.
\end{eqnarray*}
Last equalities are due to Lemma 8 from Butucea and Tsybakov \cite{ButuTsyba04}. We note that, the variance term with $\delta_{opt}$ is bigger than the variance term with $\delta^*$  but these terms are asymptotically negligible w.r.t. the bias ones. This improvement does not appear in the main term of the asymptotics. Then we conclude $\varphi_n^2=L^2\exp{\left(\frac{-4\alpha}{\delta_{opt}^r}\right)}(1+o(1)).$ The lower bound is proven in last section.
%
\paragraph{\textbf{Sketch of proof of upper bound in Theorem~\ref{theo:1}-\eqref{vitesse1}}}
For $r=2$ and $\frac{1-\eta}{2\eta}-2\alpha>0$, we have by \eqref{propa} and \eqref{propbdeux}: $$E[|d_n^2-d^2|^2]\leq\frac{8\eta^2/(1-\eta)^2}{\pi^2 n^2}e^{\frac{1-\eta}{\eta}\frac{1}{\delta^2}}+\frac{8L}{n\pi}\frac{\eta}{1-\eta-4\alpha\eta}e^{(\frac{1-\eta}{2\eta}-2\alpha)\frac{1}{\delta^2}} +L^2e^{-4\alpha/\delta^2}.$$ To select the bandwidth, we choose $\delta=\delta^*$ as solution of
$$\delta^*=\arg\inf_{\delta>0}\left\{\frac{C_{V1}}{n^2}e^{\frac{1-\eta}{\eta}\frac{1}{\delta^2}}+\frac{C_{V2}}{n}e^{(\frac{1-\eta}{2\eta}-2\alpha)\frac{1}{\delta^2}}+C_Be^{-4\alpha/\delta^2}\right\}$$ By taking derivatives, we found $\delta^*$, a positive real number satisfying $\frac{1}{\delta^{*2}}=\frac{\log n}{\frac{1-\eta}{2\eta}+2\alpha}$ we get the rate $\varphi_n^2=n^{\frac{-4\alpha}{\frac{1-\eta}{2\eta}+2\alpha}}.$ The proof of the lower bound is in last section.
%
%
\paragraph{\textbf{Proof of the parametric rate in Theorem~\ref{theo:2}}}
For $r=2$ and $\frac{1-\eta}{2\eta}-2\alpha<0$ we have by \eqref{propa} and \eqref{propbtrois}:
\begin{eqnarray*}
E[|d_n^2-d^2|^2]&\leq&\frac{8\eta^2/(1-\eta)^2}{\pi^2 n^2}e^{\frac{1-\eta}{\eta}\frac{1}{\delta^2}}+\frac{8\eta L}{4\alpha\eta-1+\eta}\cdot\frac{1}{n}+L^2e^{-4\alpha/\delta^2}.
\end{eqnarray*}
And we can write by taking $\frac{1}{\delta^{*2}}=\frac{\eta\log n}{1-\eta}$
\begin{eqnarray*}
\sup_{W_\rho\in\mathcal{A}(\alpha,2,L)}E[|d_n^2-d^2|^2]&\leq& C_V\frac{e^{\frac{1-\eta}{\eta}\frac{1}{\delta^2}}}{n^2}+C_Be^{-4\alpha/\delta^2} \leq C_V\frac{1}{n}+C_Bn^{\frac{-4\alpha}{(1-\eta)/\eta}}\\
&\leq &C_V\frac{1}{n}(1+o(1)).
\end{eqnarray*}
So we find a parametric rate. The proof of the asymptotic normality is in the section~\ref{sec:8.3}.

\paragraph{\textbf{Proof of upper bound in Theorem~\ref{theo:3}}}
Our proof is based on results of Butucea and Tsybakov \cite{ButuTsyba04}. Define $ a:=\frac{1-\eta}{4\eta}$.
\paragraph{Over the set $\Theta_1$\\}
As $0<r/2<1/2$  it is easy to remark $-(\frac{\log n}{2a}-\sqrt{\frac{\log n}{2a}})^{r/2}>-\frac{a}{2\alpha}\sqrt{\frac{\log n}{2a}}$ for $n$ large enough, and thus $\exp\left(-\frac{4\alpha}{(\delta^1_{ad})^r}\right)\geq\exp\left(-2a\sqrt{\frac{\log n}{2a}}\right).$ On the other hand the first and second variance terms found in \eqref{propbun} are equal respectively
\begin{eqnarray*}
\frac{1}{n}\exp\left(\frac{2a}{(\delta^1_{ad})^2}-\frac{2\alpha}{(\delta^1_{ad})^r}\right)&=&\exp\left(-2a\sqrt{\frac{\log n}{2a}}-2\alpha\left(\frac{\log n}{2a}-\sqrt{\frac{\log n}{2a}}\right)^{r/2}\right)\\
\frac{1}{n^2}\exp\left(\frac{4a}{(\delta^1_{ad})^2}\right)&=&\exp\left(-4a\sqrt{\frac{\log n}{2a}}\right).
\end{eqnarray*}
Therefore, with the bandwidth $\delta^1_{ad}$ the ratio of the bias term found in \eqref{propa} to the first and second variance terms are bounded from below respectively by $\exp\left(2\alpha\left(\frac{\log n}{2a}-\sqrt{\frac{\log n}{2a}}\right)^{r/2}\right)$ and $\exp\left(2a\sqrt{\frac{\log n}{2a}}\right).$ These expressions tend to $\infty$ as $n\rightarrow\infty$. Thus, the variance terms are asymptotically negligible w.r.t the bias term.It remains to check that the bias term with the bandwidth $\delta^1_{ad}$ is asymptotically bounded by $\varphi_n^2$. For $n$ large enough
\begin{eqnarray*}
L^2\exp\left(-\frac{4\alpha}{(\delta^1_{ad})^r}\right)&=&L^2\exp\left(-4\alpha\left(\frac{\log n}{2a}-\sqrt{\frac{\log n}{2a}}\right)^{r/2}\right)\\
&=&L^2\exp\left(-4\alpha(\frac{\log n}{2a})^{r/2}\left(1-(\frac{\log n}{2a})^{-1/2}\right)^{r/2}\right)\\
&\leq&L^2\exp\left(-4\alpha(\frac{\log n}{2a})^{r/2}+c(\frac{\log n}{2a})^{r/2-1/2}\right) \leq\varphi_n^2(1+o(1)).
\end{eqnarray*}
\paragraph{Over the set $\Theta_2$\\} As $r=1$ a simple calculation shows that $\delta_{opt}=\left(\frac{\log n}{2a}-\frac{\alpha}{a}\sqrt{\frac{\log n}{2a}}\right)^{-1/2}$ is a correct approximation in this case, giving a variance infinitely smaller than the bias which is of order
\begin{eqnarray*}
\varphi_n^2&=&L^2\exp\left(-\frac{4\alpha}{(\delta_{opt})}\right)=L^2\exp\left(-4\alpha\left(\frac{\log n}{2a}-\frac{\alpha}{a}\sqrt{\frac{\log n}{2a}}\right)^{1/2}\right)\\
&=&L^2\exp\left(-4\alpha\sqrt{\frac{\log n}{2a}}+\frac{2\alpha^2}{a}\right)(1+o(1)).
\end{eqnarray*}
As for the estimator with bandwidth $\delta^2_{ad}$ we get
\begin{eqnarray*}
L^2\exp\left(-\frac{4\alpha}{(\delta^2_{ad})}\right)&=&L^2\exp\left(-4\alpha\sqrt{\frac{\log n}{2a}}+\frac{2A\alpha}{a}\right)(1+o(1))\\
&=& C_2L^2\exp\left(-4\alpha\sqrt{\frac{\log n}{2a}}+\frac{2\alpha^2}{a}\right)(1+o(1)).
\end{eqnarray*}
Hence the results.

\section{Proof of the Proposition 1}
\label{sec:8}

\paragraph{Main tools}
Note that extensive use is made of formulaes \eqref{fourierp}, \eqref{fourierproun}, \eqref{fourierprodeux}, \eqref{fourierbruit} and \eqref{ens fctnel} and that Plancherel formula writes,
\begin{eqnarray}
\label{distance}
d^2:=\int_{\mathbb{R}^2}W_\rho^2(p,q)dpdq 
=\frac{1}{(2\pi)^{2}}\int_{\mathbb{R}^2} |\widetilde{W}_\rho(u,v)|^2dudv.
\end{eqnarray}

\subsection{Proof of Proposition 1-\eqref{propa}}
\label{sec:8.1}

We write $E[\cdot]$ instead of $E_\rho[\cdot]$. Because $Y_k$ and $Y_\ell$ are i.i.d for all $k\neq \ell$
\begin{eqnarray*}
E[d_n^{^2}]&=&\frac{1}{(2\pi)^{2}}\frac{1}{n(n-1)}\sum_{k\neq\ell=1}^n\int_{|t|\leq\frac{1}{\delta\sqrt{\eta}}}\int_0^\pi
\eta|t|E\left[e^{\frac{t^2}{2}(1-\eta)}e^{itY_k-itY_\ell}\right]d\phi dt\\
&=&\frac{1}{(2\pi)^{2}}\int_{|t|\leq\frac{1}{\delta\sqrt{\eta}}}\int_0^\pi\eta|t|e^{\frac{t^2}{2}(1-\eta)}|E\left[e^{itY}\right]|^2d\phi dt\\
&=&\frac{1}{(2\pi)^{2}}\int_{|t|\leq\frac{1}{\delta\sqrt{\eta}}}\int_0^\pi\eta|t|e^{\frac{t^2}{2}(1-\eta)}|\mathcal{F}[p^{\eta}_\rho(\cdot/\phi)](t)|^2d\phi dt.
\end{eqnarray*}
Use \eqref{fourierprodeux}, and a changes of variables $T=t\sqrt{\eta}$ and next the polar coordinates $u=T\cos \phi$, $v=T\sin\phi$
\begin{eqnarray*}
E[d_n^{^2}] 
&=&\frac{1}{(2\pi)^{2}}\int_{|T|\leq\frac{1}{\delta}}\int_0^\pi\eta|T||\mathcal{F}[p_\rho(\cdot/\phi)](T)|^2d\phi dT\\
\label{esperance}
&=&\frac{1}{(2\pi)^{2}}\int\int_{\|(u,v)\|^2_2\leq\frac{1}{\delta}}|\widetilde{W}_\rho(u,v)|^2dudv.
\end{eqnarray*}
So by combining \eqref{distance} et \eqref{esperance} and define $w:=(u,v)$
\begin{eqnarray*}
|E[d_n^{^2}]-d^{2}|&=&|\frac{1}{(2\pi)^{2}}\int_{|t|>\frac{1}{\delta\sqrt{\eta}}}\int_0^\pi\eta|t|e^{\frac{t^2}{2}(1-\eta)}|E[e^{itY}]|^2d\phi dt|\\
&=&\frac{1}{(2\pi)^{2}}|\int_{\mathbb{R}^2}\mathbb{I}_{\|w\|_2>1/\delta}|\widetilde{W}_\rho(w)|^2dw|\\
&\leq& \frac{1}{(2\pi)^{2}}e^{-2\alpha/\delta^r}\int_{\mathbb{R}^2}|\widetilde{W}_\rho(w)|^2e^{2\alpha\|w\|_2^r}dw\leq 
L e^{-2\alpha/\delta^r},
\end{eqnarray*}
as $W_\rho$ belongs to $\mathcal{A}(\alpha,r,L)$.

\subsection{Proof of Proposition 1-\eqref{propbun}-\eqref{propbdeux}-\eqref{propbtrois}}
\label{sec:8.2}

First we center variables
\begin{eqnarray*}
d_n^{2}-E[d_n^{2}]&=&\frac{1}{4\pi^{2}n(n-1)}\sum_{k\neq\ell=1}^n\int_{|t|\leq\frac{1}{\delta\sqrt{\eta}}}\int_0^\pi
\eta|t|e^{\frac{t^2}{2}(1-\eta)}\left(e^{itY_k-itY_\ell}-E[e^{itY}]E[e^{-itY}]\right)d\phi dt\\
&=&\frac{1}{4\pi^{2}n(n-1)}\sum_{k\neq\ell=1}^n\int_{|t|\leq\frac{1}{\delta\sqrt{\eta}}}\int_0^\pi\eta|t|e^{\frac{t^2}{2}(1-\eta)}\left(e^{itY_k}-E[e^{itY}]\right)\\
&&\cdot\left(e^{-itY_\ell}-E[e^{-itY}]\right)d\phi dt+\frac{1}{4\pi^{2}n}\sum_{k=1}^n\int_{|t|\leq\frac{1}{\delta\sqrt{\eta}}}\int_0^\pi
\eta|t|e^{\frac{t^2}{2}(1-\eta)}\left(e^{itY_k}E[e^{-itY}]\right.\\
&&\left.+e^{-itY_k}E[e^{itY}]\right)d\phi dt-2|E[e^{itY}]|^2.
\end{eqnarray*}
Let define by $Z_k(t)=Z_k:=e^{itY_k}-E[e^{itY}]$, and $\bar{Z}_k$ its complex conjugate, then:
\begin{eqnarray*}
d_n^{2}-E[d_n^{2}]&=&\frac{1}{(2\pi)^{2}}\left(\frac{1}{n(n-1)}\sum_{k\neq\ell}\int_{|t|\leq\frac{1}{\delta\sqrt{\eta}}}\int_0^\pi \eta|t|e^{\frac{t^2}{2}(1-\eta)}Z_k\bar{Z}_\ell d\phi dt\right.\\
&&+\left.\frac{1}{n}\sum_{j}\int_{|t|\leq\frac{1}{\delta\sqrt{\eta}}}\int_0^\pi\eta|t|e^{\frac{t^2}{2}(1-\eta)}{\left(Z_jE[e^{-itY}]+\bar{Z}_jE[e^{itY}]\right)}d\phi dt\right).
\end{eqnarray*}
Denote by $J_1$ and $J_2$ respectively the first and the second term of the previous sum, we have then
\begin{equation}
\label{variance}
Var(d_n^{2})=E[(d_n^{2}-E[d_n^{2}])^2]=E[J_1^2]+E[J_2^2]+2E[J_1J_2].
\end{equation}
See that the third part of the  previous sum:
\begin{eqnarray*}
E[J_1J_2]&=&\frac{1}{(2\pi)^{4}}\frac{1}{n^2(n-1)}\sum_{k\neq\ell}\sum_{j}E\left[\left(\int_{|t|\leq\frac{1}{\delta\sqrt{\eta}}}\int_0^\pi \eta|t|e^{\frac{t^2}{2}(1-\eta)}Z_k\bar{Z}_{\ell}d\phi dt\right)\right.\\
&&\cdot\left.\left(\int_{|t|\leq\frac{1}{\delta\sqrt{\eta}}}\int_0^\pi\eta|t|e^{\frac{t^2}{2}(1-\eta)}\left(E[e^{-itY}]Z_j
+E[e^{itY}]\bar{Z}_{j}\right)d\phi dt\right)\right]=0.
\end{eqnarray*}
By noticing $E[Z_j]=0$ for all $j=1,...,n$, and because there always exists a $j\neq k$ and thus $Z_j$, $Z_k$ are independent  or a $j\neq \ell$ and $Z_j$, $Z_\ell$ are independent. Now study
\begin{eqnarray*}
E[J_1^2]&=&\frac{1}{16\pi^{4}n^2(n-1)^2}E\left[\left(\sum_{k\neq\ell}\int_{|t|\leq\frac{1}{\delta\sqrt{\eta}}}\int_0^\pi
\eta|t|e^{\frac{t^2}{2}(1-\eta)}Z_{k}\bar{Z}_{\ell}d\phi dt\right)^2\right]\\
&=&\frac{1}{16\pi^{4}n^2(n-1)^2}\sum_{k_1\neq\ell_1}\sum_{k_2\neq\ell_2}E\left[\left(\int_{|t|\leq\frac{1}{\delta\sqrt{\eta}}}\int_0^\pi\eta|t|e^{\frac{t^2}{2}(1-\eta)}Z_{k_1}\bar{Z}_{\ell_1}d\phi 
dt\right)\right.\\
&&\cdot\left.\left(\int_{|t|\leq\frac{1}{\delta\sqrt{\eta}}}\int_0^\pi \eta|t|e^{\frac{t^2}{2}(1-\eta)}Z_{k_2}\bar{Z}_{\ell_2}d\phi dt\right)\right].
\end{eqnarray*}
Note that,as soon as an indices $k_1$ ,$\ell_1$, $k_2$ ,$\ell_2$ is different from the others, the expected value is $0$. Thus,
\begin{eqnarray*}
E[J_1^2]&=&\frac{1}{16\pi^{4}n^2(n-1)^2}\left(\sum_{k\neq\ell}E\left[\left(\int_{|t|\leq\frac{1}{\delta\sqrt{\eta}}}\int_0^\pi\eta|t| e^{\frac{t^2}{2}(1-\eta)}Z_{k}\bar{Z}_{\ell}d\phi 
dt\right)^2\right]\right.\\
&&\left.+\sum_{k\neq\ell}E\left[|\int_{|t|\leq\frac{1}{\delta\sqrt{\eta}}}\int_0^\pi\eta|t|e^{\frac{t^2}{2}(1-\eta)}Z_{k}\bar{Z}_{\ell}d\phi dt|^2\right]\right)
\end{eqnarray*}
\begin{eqnarray*}
E[J_1^2]&=&\frac{1}{16\pi^{4}n(n-1)}\cdot\frac{1}{2}E\left[\left(\int_{|t|\leq\frac{1}{\delta\sqrt{\eta}}}\int_0^\pi\eta|t|e^{\frac{t^2}{2}(1-\eta)}Z_1\bar{Z}_{2}d\phi dt\right)^2\right.\\
&&\left.+\left(\int_{|t|\leq\frac{1}{\delta\sqrt{\eta}}}\int_0^\pi\eta|t|e^{\frac{t^2}{2}(1-\eta)}Z_2\bar{Z}_{1}d\phi dt\right)^2\right]\\
&&+\frac{1}{16\pi^{4}n(n-1)}E\left[|\int_{|t|\leq\frac{1}{\delta\sqrt{\eta}}}\int_0^\pi\eta|t|e^{\frac{t^2}{2}(1-\eta)}Z_1\bar{Z}_{2}d\phi dt|^2\right]
\end{eqnarray*}
\begin{eqnarray*}
E[J_1^2]&=&\frac{1}{16\pi^{4}n(n-1)}E\left[\Re e\left(\left(\int_{|t|\leq\frac{1}{\delta\sqrt{\eta}}}\int_0^\pi\eta|t|e^{\frac{t^2}{2}(1-\eta)}Z_1\bar{Z}_{2}d\phi 
dt\right)^2\right)\right]\\
&&+\frac{1}{16\pi^{4}n(n-1)}E\left[|\int_{|t|\leq\frac{1}{\delta\sqrt{\eta}}}\int_0^\pi\eta|t|e^{\frac{t^2}{2}(1-\eta)}Z_1\bar{Z}_{2}d\phi dt|^2\right].
\end{eqnarray*}
By noticing that $|\Re e(z)|\leq |z|$  and using the fact $|Z_k|\leq 2$, we get
\begin{eqnarray}
E[J_1^2]&\leq&\frac{1}{8\pi^{4}n(n-1)}E\left[|\int_{|t|\leq\frac{1}{\delta\sqrt{\eta}}}\int_0^\pi\eta|t|
e^{\frac{t^2}{2}(1-\eta)}Z_1\bar{Z}_{2}d\phi dt|^2\right]\nonumber\\
&\leq&\frac{2\pi^2}{\pi^{4}n^2}\left(\int_{|t|\leq\frac{1}{\delta\sqrt{\eta}}}\eta|t|e^{\frac{t^2}{2}(1-\eta)}d\phi dt\right)^2\nonumber \\
\label{Jun}
&\leq &\frac{8\eta^2}{\pi^{2}(1-\eta)^2}\cdot\frac{1}{n^2}\exp\left(\frac{1-\eta}{\eta\delta^2}\right).
\end{eqnarray}
Noticing that $(Z_k)_k$ are i.i.d and centered, we then have by developing the square:
\begin{eqnarray}
E[J_2^2]&=&E\left[\frac{1}{(4\pi^2 n)^{2}}\left(\sum_{k=1}^n\int_{|t|\leq\frac{1}{\delta\sqrt{\eta}}}\int_0^\pi
\eta|t|e^{\frac{t^2}{2}(1-\eta)}2\Re e(E[e^{itY}]\bar{Z}_k)d\phi dt\right)^2\right]\nonumber\\
&=&\frac{1}{16\pi^{4}n}E\left[\left(\int_{|t|\leq\frac{1}{\delta\sqrt{\eta}}}\int_0^\pi\eta|t|e^{\frac{t^2}{2}(1-\eta)}2\Re 
e(E[e^{itY}]\bar{Z}_1)d\phi dt\right)^2\right] \nonumber\\
\label{termdom}
&=&\frac{1}{4\pi^{4}n}E\left[\left(\int_{|t|\leq\frac{1}{\delta\sqrt{\eta}}}\int_0^\pi\eta|t|e^{\frac{t^2}{2}(1-\eta)}\Re 
e(E[e^{itY}]\bar{Z}_1)d\phi dt\right)^2\right]\\
&\leq&\frac{1}{\pi^{4}n}\left(\int_{|t|\leq\frac{1}{\delta\sqrt{\eta}}}\int_0^\pi\eta|t|e^{\frac{t^2}{2}(1-\eta)}|E[e^{itY}]|d\phi  dt\right)^2, \nonumber
\end{eqnarray}
as $|\Re e(z)|\leq |z|$ and $|\bar{Z}_1|\leq 2$. Then use successively \eqref{fourierproun}and \eqref{fourierbruit}, next \eqref{fourierprodeux},and a change of variables $T=t\sqrt{\eta}$
\begin{eqnarray*}
E[J_2^2]&\leq&\frac{1}{\pi^{4}n}\left(\int_{|t|\leq\frac{1}{\delta\sqrt{\eta}}}\int_0^\pi\eta|t|e^{\frac{t^2}{2}(1-\eta)}|E[e^{itY}]|d\phi dt\right)^2\\
&=&\frac{1}{\pi^{4}}\frac{1}{n}\left(\int_{|t|\leq\frac{1}{\delta\sqrt{\eta}}}\int_0^\pi\eta|t|e^{\frac{t^2}{4}(1-\eta)}\frac{|\textit{F}_1[p_\rho^\eta(./\phi)](t)|}{|\widetilde{N}^\eta(t)|}d\phi dt\right)^2\\
&=&\frac{1}{\pi^{4}}\frac{1}{n}\left(\int_{|t|\leq\frac{1}{\delta\sqrt{\eta}}}\int_0^\pi\eta|t|e^{\frac{t^2}{4}(1-\eta)}|\textit{F}_1[p_\rho(./\phi)](\sqrt\eta t)|d\phi dt\right)^2\\
&=&\frac{1}{\pi^{4}}\frac{1}{n}\left(\int_0^\pi \int_{|T|<\frac{1}{\delta}}|T|e^{\frac{T^2}{4\eta}(1-\eta)}|\textit{F}_1[p_\rho(./\phi)](T)|d\phi dT\right)^2.
\end{eqnarray*}
Then, by \eqref{fourierp} and next use the polar coordinates $u=T\cos\phi$, $v=T\sin\phi$
\begin{eqnarray*}
E[J_2^2]&=&\frac{1}{\pi^{4}}\frac{1}{n}\left(\int_0^\pi\int_{|T|<\frac{1}{\delta}}|T|e^{\frac{T^2}{4\eta}(1-\eta)}|\widetilde{W}_{\rho}(T\cos\phi,T\sin\phi)|d\phi dT\right)^2\\
&=&\frac{1}{\pi^{4}}\frac{1}{n}\left(\int_{\|(u,v)\|_2\leq 1/\delta}e^{\frac{1-\eta}{4\eta}\|(u,v)\|^2_2}|\widetilde{W}_{\rho}(u,v)|dudv\right)^2.
\end{eqnarray*}
Define $z:=(u,v)$, and use Cauchy-Schwartz inequality and \eqref{ens fctnel}
\begin{eqnarray}
E[J_2^2]&=&\frac{1}{\pi^{4}}\frac{1}{n}\left(\int_{\|z\|_2\leq 1/\delta}e^{\frac{1-\eta}{4\eta}\|z\|^2_2}|\widetilde{W}_{\rho}(z)|dz\right)^2\nonumber\\
&\leq&\frac{1}{n\pi^4}(2\pi)^2L\int_{\|z\|_2\leq1/\delta}te^{\frac{1-\eta}{2\eta}\|z\|_2^2-2\alpha\|z\|_2^r}dz\nonumber\\
\label{Jdeux} &\leq&\frac{8L}{n\pi} \int_0^{1/\delta}te^{\frac{1-\eta}{2\eta}t^2-2\alpha t^r}dt.
\end{eqnarray}
\begin{enumerate}
\item For $0<r<2$ and according to Lemma 6 of Butucea and Tsybakov \cite{ButuTsyba04} we get:
\begin{eqnarray}
\label{run}
\frac{8L}{n\pi}\int_0^{1/\delta}te^{\frac{1-\eta}{2\eta}t^2-2\alpha t^r}dt\leq\frac{8L}{n\pi}\frac{\eta}{1-\eta}e^{\frac{1-\eta}{2\eta}\frac{1}{\delta^2}-2\alpha\frac{1}{\delta^r}}.
\end{eqnarray}
The expressions \eqref{Jun} and \eqref{Jdeux} together with \eqref{run} conclude \eqref{propbun}.\item For $r=2$ and $\frac{1-\eta}{2\eta}-2\alpha>0$ and according to Lemma 6 of \cite{ButuTsyba04} we get:
\begin{eqnarray}
\frac{8L}{n\pi}\int_0^{1/\delta}te^{\frac{1-\eta}{2\eta}t^2-2\alpha t^2}dt&\leq&\frac{8L}{n\pi}\cdot\frac{1}{2(\frac{1-\eta}{2\eta}-2\alpha)}e^{(\frac{1-\eta}{2\eta}-2\alpha)\frac{1}{\delta^2}}\nonumber\\
\label{rdeux}
&\leq&\frac{8L}{n\pi}\cdot\frac{\eta}{1-\eta-4\alpha\eta}e^{(\frac{1-\eta}{2\eta}-2\alpha)\frac{1}{\delta^2}}.
\end{eqnarray}
The expressions \eqref{Jun} and \eqref{Jdeux} together with \eqref{rdeux} conclude \eqref{propbdeux}.\item For $r=2$ and $\frac{1-\eta}{2\eta}-2\alpha<0$ we have:
\begin{eqnarray}
\frac{8L}{n\pi}\int_0^{1/\delta}te^{\frac{1-\eta}{2\eta}t^2-2\alpha          t^2}dt&\leq&\frac{4L}{2\alpha-\frac{1-\eta}{2\eta}}\frac{1}{n}
\label{rtrois} \leq\frac{8\eta L}{4\alpha\eta-1+\eta}\cdot\frac{1}{n}.
\end{eqnarray}
The expressions \eqref{Jun} and \eqref{Jdeux} together with \eqref{rtrois} conclude \eqref{propbtrois}.
\end{enumerate}

\subsection{the asymptotic normality}
\label{sec:8.3}

Let $r=2$, $\frac{1-\eta}{2\eta}-2\alpha<0$ and $\delta=\delta^*=\left(\frac{\eta\log n}{1-\eta}\right)^{-1/2}$.
$$\sqrt{n}(d_n^2-d^2)=\sqrt{n}(d_n^2-E[d_n^2])+\sqrt{nB(d_n^2)}.$$
The term $\sqrt{nB(d_n^2)}\leq\sqrt{n}Le^{-\frac{2\alpha}{(\delta^*)^r}}$ tends to 0 as $n\rightarrow \infty$.\\
Moreover $\sqrt{n}(d_n^2-E[d_n^2])=\sqrt{n}(J_1+J_2)$, where $J_{1,2}$ are centered and were defined in Section~\ref{sec:8.2}. It has been shown in \eqref{rtrois} that the dominating term in the variance $E[(d_n^{2}-E[d_n^{2}])^2]$ is given by $E[J^2_2]$ defined in \eqref{termdom}. That means $nE[J^2_1]=o(1)$, as $n\rightarrow \infty$ and $\sqrt{n}J_1\stackrel{P}{\rightarrow} 0$. Thus, the asymptotic normality is given by the term $\sqrt{n}J_2$.\\
As $J_2=\left.\frac{1}{n}\sum_{j}\int_{|t|\leq\frac{1}{\delta\sqrt{\eta}}}\int_0^\pi\eta|t|e^{\frac{t^2}{2}(1-\eta)}{\left(Z_jE[e^{-itY}]+\bar{Z}_jE[e^{itY}]\right)}d\phi dt\right)$, we can use a classical central limit theorem for i.i.d. random 
variables with finite variance and the asymptotic variance is given by the limit of $n E[J_2^2]$. Let us study $\lim_{n\to\infty} nE[J^2_2]$
\begin{eqnarray*}
n E[J_2^2]&=&\frac{1}{4\pi^{4}}E\left[\Re e\left(\int_{|t|\leq\frac{1}{\delta\sqrt{\eta}}}\int_0^\pi \eta|t|e^{\frac{t^2}{2}(1-\eta)}(E[e^{itY}]e^{-itY}\right.\right.\\
&&\left.\left.-E[e^{itY}]E[e^{-itY}])d\phi dt\right)^2\right]\\
&=&\frac{1}{4\pi^{4}}E\left[\Re e\left(\int_{|t|\leq\frac{1}{\delta\sqrt{\eta}}}\int_0^\pi\eta|t|e^{\frac{t^2}{2}(1-\eta)}E[e^{itY}]e^{-itY}d\phi
dt\right)^2\right]\\
&&-\frac{1}{4\pi^{4}}\left(\int_{|t|\leq\frac{1}{\delta\sqrt{\eta}}}\int_0^\pi\eta|t|e^{\frac{t^2}{2}(1-\eta)}E[e^{itY}]E[e^{-itY}]d\phi dt\right)^2:=A_1-A_2.
\end{eqnarray*}
On the one hand, we already proved in section~\ref{sec:8.1}, $A_2=4\left(E[d_n^2]\right)^2$. Therefore, $\lim_{n\to\infty}A_2=\left\|W_\rho\right\|^2_2$. On the other hand
\begin{eqnarray*}
A_1&=&\frac{1}{4\pi^{4}}E\left[\Re 
e\left(\int_{|t|\leq\frac{1}{\delta\sqrt{\eta}}}\int_0^\pi\eta|t|e^{\frac{t^2}{2}(1-\eta)}E[e^{itY}]e^{-itY}d\phi
dt\right)^2\right]\\
&=&\frac{1}{4\pi^{4}}\int_{|t_1|\leq\frac{1}{\delta\sqrt{\eta}}}\int_{|t_2|\leq\frac{1}{\delta\sqrt{\eta}}}\int_0^\pi\int_0^\pi\eta^2|t_1||t_2|e^{\frac{t_1^2+t_2^2}{2}(1-\eta)}E[e^{it_1Y}]E[e^{it_2Y}]\\
&&\cdot E[e^{-i(t_1+t_2)Y}]d\phi_1 dt_1d\phi_2 dt_2.
\end{eqnarray*}
By  changing the variable $t$ into $t/\sqrt{\eta}$ and as $Y/\sqrt{\eta}=(X+\sqrt{\frac{1-\eta}{2\eta}}\xi)$ we get
\begin{eqnarray*}
A_1&=&\frac{1}{4\pi^{2}}\int_{|t_1|\leq\frac{1}{\delta}}\int_{|t_2|\leq\frac{1}{\delta}}|t_1||t_2|e^{\frac{1-\eta}{2\eta}(t_1^2+t_2^2)}E[e^{it_1Y/\sqrt{\eta}}]E[e^{it_2Y/\sqrt{\eta}}]\\
&&\cdot E[e^{-i(t_1+t_2)Y/\sqrt{\eta}}]dt_1dt_2\\
&=&\frac{1}{4\pi^{2}}\int_{|t_1|\leq\frac{1}{\delta}}\int_{|t_2|\leq\frac{1}{\delta}}|t_1||t_2|e^{\frac{1-\eta}{2\eta}(t_1^2+t_2^2)}E[e^{it_1(X+\sqrt{\frac{1-\eta}{2\eta}}\xi)}]E[e^{it_2(X+\sqrt{\frac{1-\eta}{2\eta}}\xi)}]\\
&&\cdot E[e^{-i(t_1+t_2)(X+\sqrt{\frac{1-\eta}{2\eta}}\xi)}]dt_1 dt_2.
\end{eqnarray*}
As $X$ and $\xi$ are independent and since $E[e^{iT\sqrt{\frac{1-\eta}{2\eta}}\xi)}]=e^{-T^2\frac{1-\eta}{4\eta}}$, we get
\begin{eqnarray*}
A_1&=&\frac{1}{4\pi^{2}}\int_{|t_1|\leq\frac{1}{\delta}}\int_{|t_2|\leq\frac{1}{\delta}}|t_1||t_2|e^{\frac{1-\eta}{2\eta}(t_1^2+t_2^2)}e^{-\frac{1-\eta}{4\eta}(t_1^2+t_2^2)}e^{-\frac{1-\eta}{4\eta}(t_1+t_2)^2}E[e^{it_1X}]\\
&&\cdot E[e^{it_2X}] E[e^{-i(t_1+t_2)X}]dt_1 dt_2\\
&=&\frac{1}{4\pi^{2}}\int_{|t_1|\leq\frac{1}{\delta}}\int_{|t_2|\leq\frac{1}{\delta}}|t_1||t_2|e^{\frac{1-\eta}{2\eta}t_1t_2}E[e^{it_1X}]E[e^{it_2X}] E[e^{-i(t_1+t_2)X}]dt_1 dt_2,
\end{eqnarray*}
and 
$\lim_{n\to\infty}A_1=\frac{1}{4\pi^{2}}\int\int|t_1||t_2|e^{\frac{1-\eta}{2\eta}t_1t_2}E[e^{it_1X}]E[e^{it_2X}] E[e^{-i(t_1+t_2)X}]dt_1 dt_2$.\\
By denoting $\mathcal{W}=\lim_{n\to\infty}(A_1-A_2)$, we get the result.

\section{Proofs of lower bounds}
\label{sec:9}

In this section, we will show the lower bounds of Theorem~\ref{theo:1}. For that we will be based on the results of Butucea and Tsybakov \cite{ButuTsyba04} . They show that the problem of bound from above the minimax risk can be reduce to two functions $W_{\rho_1}$ and $W_{\rho_0}$ depending on a parameter $\tilde{\delta_n}=\tilde{\delta}$ such that $\tilde{\delta}\rightarrow 0$ as $n\rightarrow 0$. The choice of $\tilde{\delta}$ insures the existence of the lower bound. The parameter $\tilde{\delta}$ is the unique solution of the equation
\begin{eqnarray}
\label{equah}
\frac{2\alpha}{\tilde{\delta}^r}+\frac{1-\eta}{2\eta\tilde{\delta}^2}=\log n+(\log\log n)^2.
\end{eqnarray}
If $0<r<2$, notice that it is different of the $\delta$ appearing in the expression of our estimator defined in \eqref{estimateur}. And for $r=2$, we take
\begin{eqnarray}
\label{equahdeux} \tilde{\delta}=\left(\frac{\log (n\log n)}{2(a+\alpha)}\right)^{-1/2},\quad \text{where}\quad a=\frac{1-\eta}{4\eta}.
\end{eqnarray}
We will use Wigner functions $W_{\rho_0}$ and $W_{\rho_1}$ built by Butucea \textit{et al.} (2007) \cite{ButGutArt05} in their first prepublication like certain results coming from this construction. $W_{\rho_0}$ is a fixed function corresponding to the density matrix $\rho_o$, and $W_{\rho_1}$ is of the form
$$W_{\rho_1}(z)=W_{\rho_0}(z)+V_{\tilde{\delta}}(z)\quad \text{and}$$
$$\rho_1=\rho_0 +\tau^{\tilde{\delta}}$$
such that $\rho_1$ is a density matrix (positive and trace equal to one) with Radon transforms $p_1$. Note that the function $V_{\tilde{\delta}}$ is not a Wigner function of a density matrix but belongs to the linear span of the space of Wigner functions and its corresponding matrix $\tau^{\tilde{\delta}}$ is in the linear span of density matrix.  We will 
detail in a next paragraph the construction of $W_{\rho_{0,1}}$, $\rho_{0,1}$ and $V_{\tilde{\delta}}$ as well as the results which results from this. As we have stipulated it higher, we will use lemma 4 in Butucea and Tsybakov \cite{ButuTsyba04}. Let us suppose first of all that the following conditions are satisfied:
\begin{equation}
\label{hyp1} W_{\rho_1},W_{\rho_0}\in\mathcal{A}(\alpha,r,L),
\end{equation}
\begin{equation}
\label{hyp2} |d_1^2-d_0^2|=|\|W_{\rho_1}\|^2_2-\|W_{\rho_0}\|^2_2|\geq 2\phi_n(1+o(1)),\;n\rightarrow\infty,
\end{equation}
\begin{equation}
\label{hyp3} n\chi^2:=n\int_0^\pi \int\frac{(p_1^\eta(y)-p_0^\eta(y))^2}{p_0^\eta(y)}dyd\phi=o(1),\;n\rightarrow\infty.
\end{equation}
Then we reduce the minimax risk to these two functions, $W_{\rho_1}$, $W_{\rho_0}$, and note $\widehat{d}_n^2$ an arbitrary estimator of $d_\rho^2:=\|W_\rho\|^2_2$, then we get for some $0<\tau<1$
\begin{eqnarray*}
\inf_{\widehat{d}_n^2}\sup_{W_\rho\in\mathcal{A}(\alpha,r,L)}E[|\widehat{d}_n^2-d_\rho^2|^2]
&\geq&\inf_{\widehat{d}_n^2}\frac{1}{2}(E_{\rho_0}[|\widehat{d}_n^2-d_{\rho_0}^2|^2]+E_{\rho_1}[|\widehat{d}_n^2-d_{\rho_1}^2|^2])\\
&\geq&\inf_{\widehat{d}_n^2}\frac{1}{2}(E_{\rho_0}[|\widehat{d}_n^2-d_{\rho_0}^2|^2]\\
&&+(1-\tau)E_{\rho_0}[\mathbb{I}(\frac{dP^{\eta}_{\rho_1}}{dP^{\eta}_{\rho_0}}\geq 1-\tau)|\widehat{d}_n^2-d_{\rho_1}^2|^2])\\
&\geq&\inf_{\widehat{d}_n^2}\frac{1}{2}(1-\tau)(E_{\rho_0}[\mathbb{I}(\frac{dP^{\eta}_{\rho_1}}{dP^{\eta}_{\rho_0}}\geq
1-\tau)(|\widehat{d}_n^2-d_{\rho_0}^2|^2\\
&&+|\widehat{d}_n^2-d_{\rho_1}^2|^2)]).
\end{eqnarray*}
As $a^2+b^2\geq (a-b)^2$ for $a$ and $b$ positives reals numbers, we can get ride of the estimator.
\begin{eqnarray*}
&\geq&\frac{1}{4}(1-\tau)E_{\rho_0}[\mathbb{I}(\frac{dP^{\eta}_{\rho_1}}{dP^{\eta}_{\rho_0}}\geq(1-\tau))|d_{\rho_1}^2-d_{\rho_0}^2|^2]\\
&\geq&(1-\tau)\phi_n^2(1-P_{\rho_0}(\frac{dP^{\eta}_{\rho_1}}{dP^{\eta}_{\rho_0}}-1<-\tau))\\
&\geq&(1-\tau)\phi_n^2(1-\frac{1}{\tau^2}\int(\frac{dP^{\eta}_{\rho_1}}{dP^{\eta}_{\rho_0}}-1)^2dP^{\eta}_{\rho_0}).
\end{eqnarray*}
By supposing $n\chi^2\leq\tau^4$ the last inequality is undervalued by $(1-\tau)^2\phi_n^2(1+\tau)$. It is enough to check \eqref{hyp3}, in order to get $\tau_n\rightarrow 0$ as $n\rightarrow\infty$, and we obtain a lower bound for the minimax risk of order $\phi_n^2(1+o(1))$for any estimator $\widehat{d}_n^2$. Our proof of lower bounds is quite similar to the one of Butucea \textit{et al.} (2007) \cite{ButGutArt05}. The main difference is the proof of \eqref{hyp2} as we don't bound from below the Wigner function but the quadratic functional of the Wigner function. Nevertheless, for the reader's convenience, we reproduce key proofs to complete the proof of the lower bounds.

\subsubsection{The density matrix $\rho_0$}
\label{sec:9.0.1}

The main difference with the construction in Butucea \textit{et al.} (2007) \cite{ButGutArt05} is that they had considered two Wigner functions $W_{\rho_1}$ and $W_{\rho_2}$ with $W_{\rho_1,\rho_2}=W_{\rho_0}\pm V_{\tilde{\delta}}$ while we consider only the Wigner function $W_{\rho_1}=W_{\rho_0}+ V_{\tilde{\delta}}$. Because we have to bound from below the quantity $|\|W_{\rho_1}\|^2_2-\|W_{\rho_0}\|^2_2|$ instead of $\|W_{\rho_1}-W_{\rho_0}\|^2_2$ and we must make sure that $\widetilde{W}_{\rho_0}$ and $\widetilde{V}_{\tilde{\delta}}$ are positive functions.
In this paragraph we will recall some results and lemmas of Butucea \textit{et al.} (2007) \cite{ButGutArt05} about the density matrix $\rho_0$ and its corresponding Wigner function. They had constructed a family of density matrices $\rho^{\beta,\xi}$ from which they selected $\rho_0=\rho^{\beta_0,\xi_0}$ with Radon transform $p_{\beta}^{\xi}$ equals to
$$p_\beta^\xi(x,\phi):=\int_0^1\frac{f_\beta^\xi(z)}{\sqrt{\pi(1-z^2)}}\exp{\left(-x^2\frac{1-z}{1+z}\right)}dz,$$
where$f_\beta^\xi(z)=\beta(1-z)^\beta/(1-\xi)\mathbb{I}(\xi\leq z\leq 1)$, for some $0<\beta,\xi\leq 1$.
The Fourier transform is
$$\widetilde{W_\beta^\xi}(w)=\mathcal{F}_1[p_\beta^\xi](\|w\|,\phi)=\int_0^1\frac{f_\beta^\xi(z)}{1-z}\exp\left(-\|w\|^2\frac{1+z}{4(1-z)}\right)dz.$$
Notice that the Fourier transform is positive and $\widetilde{W_\beta^\xi}(0)=1$. The study of the asymptotic behavior of such functions is done in lemmae~\ref{lm:1} and ~\ref{lm:2}. Lemma~\ref{lm:3} proves the fact that $W_\beta^\xi$ belongs to the class $\mathcal{A}(\alpha,r,L)$ for $\beta>0$ small enough and $\xi$ close to 1.
\begin{lm}
\label{lm:1}
For all $0<\beta,\xi\leq 1$ and $|x|>1$ there exist constants $c, C$ depending on $\beta$ and $\xi$, such that $$c|x|^{-(1+2\beta)}\leq p_\beta^\xi(x,\phi)\leq C|x|^{-(1+2\beta)}.$$
\end{lm}
\begin{lm}
\label{lm:2}
For all $0<\beta,\xi\leq 1$ we have $$\rho^{\beta,\xi}_{n,n}=\frac{\beta}{(1-\xi)^\beta}\Gamma(\beta+1)n^{-(1+\beta)}(1+o(1)),\,\,n\rightarrow\infty.$$
\end{lm}
\begin{lm}
\label{lm:3}For any $(\alpha,r,L)$ such that $0<r\leq 2$, there exists an $0<\beta,\xi\leq 1$ such that $W_\beta^\xi$ belongs to the class $\mathcal{A}(\alpha,r,L)$.
\end{lm}
We refer for the proof of these lemmae to Butucea \textit{et al.} (2007) \cite{ButGutArt05}.

\subsubsection{Construction of $V_{\tilde{\delta}}$ and asymptotic properties of $\tau^{\tilde{\delta}}$}
\label{sec:9.0.2}

For using the same construction as Butucea \textit{et al.} (2007) \cite{ButGutArt05}, we have to define on $\mathbb{R}^2$ the function $V_{\tilde{\delta}}$ whose Fourier transform is
$$\mathcal{F}_2[V_{\tilde{\delta}}](w):=\widetilde{V}_{\tilde{\delta}}(w)=J_{\tilde{\delta}}(t)=2\sqrt{rL\pi\alpha}\tilde{\delta}^{(2-r)/2}e^{\alpha/\tilde{\delta}^r}e^{-2\alpha|t|^r}J(|t|^r-\frac{1}{\tilde{\delta}^r}),$$
where $t=\|w\|$, and $J$ is a 3-times continuously differentiable function on $\mathbb{R}$ with its first 3 derivatives uniformly bounded on $\mathbb{R}$ such that for any $\lambda>0$ and any $D>4\lambda$
$$\mathbb{I}(2\lambda\leq u\leq D-2\lambda)\leq J(u)\leq\mathbb{I}(\lambda\leq u\leq D-\lambda),\quad \text{for all}\,u\in\mathbb{R}.$$
We choose $\tilde{\delta}$ solution of \eqref{equah} when $0<r<2$ and $\tilde{\delta}$ such as in \eqref{equahdeux}, when $r=2$. We want $V_{\tilde{\delta}}$ to be a function of a density matrix belonging to the linear span of the space of Wigner functions and its corresponding matrix $\tau^{\tilde{\delta}}$ belonging to the linear span of density matrix. For that, we use an important property of Wigner functions: the isometry (up to a constant) between the linear span of density matrices and that of Wigner functions with respect to the $\mathbb{L}_2$-distances, in particular
$$\|W_{\rho_2}-W_{\rho_1}\|=:\int\int|W_{\rho_2}(p,q)-W_{\rho_1}(p,q)|^2dpdq=\frac{1}{2\pi}\|\rho_2-\rho_1\|^2_2,$$
for any $\rho_2$, $\rho_1$. Note that because the function $V_{\tilde{\delta}}$ is invariant under rotations in the plane, the corresponding matrix has all off-diagonal elements equal to 0 and for the diagonal ones we can use the following formula from Leonhardt (1997) \cite{Leon97}
$$\tau^{\tilde{\delta}}_{nn}=4\pi^2\int^1_0 L_n(t^2/2)e^{-t^2/4}tJ_{\tilde{\delta}}(t)dt.$$
And as our choice of $V_{\tilde{\delta}}$ is the same as Butucea \textit{et al.} (2007) \cite{ButGutArt05}, we have the same asymptotic  behavior derived in the following lemma.
\begin{lm}
\label{lm:4}
The matrix $\tau^{\tilde{\delta}}$ has the following asymptotic behavior $$\tau^{\tilde{\delta_n}}_{nn}=O(n^{-5/4})o_{\tilde{\delta}}(1).$$
\end{lm}
For the proof of this lemma we refer to Butucea \textit{et al.} (2007) \cite{ButGutArt05}. We have now to prove conditions \eqref{hyp1}, \eqref{hyp2} and \eqref{hyp3} to obtain the lower bound.

\subsection{Proof of conditions \eqref{hyp1}, \eqref{hyp2} and \eqref{hyp3}}
\label{sec:9.1}

\paragraph{\textbf{Proof of \eqref{hyp1}}}
From Lemma~ref{lm:3} we get for any $\beta$ small enough and $\xi$ close to 1 that the Wigner function $W_\beta^\xi$ belongs to the class $\mathcal{A}(\alpha,r,a^2L)$. And the Lemmae~ref{lm:2} and ~ref{lm:3} implies that for any $\beta<1/4$ the diagonal matrix $\rho_1=\rho^{\beta,\xi}+\tau^{\tilde{\delta}}$ is positive with trace one for $\tilde{\delta}$ small enough. Thus there exists an $\beta_0$, $\xi_0$ such that the corresponding matrix $\rho_1$ is a density matrix and
$W_{\rho_0}=W_{\beta_0}^{\xi_0}\in\mathcal{A}(\alpha,r,a^2L)$. Let us prove that $W_{\rho_1}\in\mathcal{A}(\alpha,r,L)$. By triangle inequality
\begin{eqnarray*}
\|\mathcal{F}_2[W_{\rho_1}]e^{\alpha\|.\|^r}\|_2 &\leq &\|\mathcal{F}_2[W_{\rho_0}]e^{\alpha\|.\|^r}\|_2+\|\mathcal{F}_2[V_{\tilde{\delta}}]e^{\alpha\|.\|^r}\|_2\\
&\leq& 2\pi a\sqrt{L}+\|\mathcal{F}_2[V_{\tilde{\delta}}]e^{\alpha\|.\|^r}\|_2.
\end{eqnarray*}
Now, by the change of variables $u=t\cos\phi$, $v=t\sin\phi$
\begin{eqnarray*}
\|\mathcal{F}_2[V_{\tilde{\delta}}]e^{\alpha\|.\|^r}\|_2^2&=&\int_{\mathbb{R}^2}|\mathcal{F}_2[V_{\tilde{\delta}}](w)|^2e^{2\alpha\|w\|^r}dw\\
&=&\int_0^\pi\int_{\mathbb{R}}|t||\mathcal{F}_2[V_{\tilde{\delta}}](t\cos\phi,t\sin\phi)|^2e^{2\alpha|t|^r}dtd\phi\\
&=&\pi\int_{\mathbb{R}}|t||\mathcal{F}_2[V_{\tilde{\delta}}](t\cos\phi,t\sin\phi)|^2e^{2\alpha|t|^r}dt\\
&=&\pi\int_{\mathbb{R}}|t||J_{\tilde{\delta}}(t)|^2e^{2\alpha|t|^r}dt\\
&\leq& 2^2\pi^2L\alpha r\tilde{\delta}^{2-r}e^{2\alpha/\tilde{\delta}^r}2\int^\infty_{(\lambda+\frac{1}{\tilde{\delta}^r})^{1/r}}
t e^{-2\alpha t^r}dt\leq 2^2\pi^2Le^{-2\alpha\lambda}.
\end{eqnarray*}
Thus, it is enough to take $a=1-e^{-\alpha\lambda/2}$ to get $W_{\rho_1}\in\mathcal{A}(\alpha,r,L(1-e^{-\alpha\lambda/2}+e^{-\alpha\lambda})^2) \subset\mathcal{A}(\alpha,r,L).$
\paragraph{\textbf{Proof of \eqref{hyp2}}}
By noticing that $\widetilde{W}_{\rho_0}$ and $\widetilde{V}_{\tilde{\delta}}$ are positive functions we get
\begin{eqnarray*}
|\|W_{\rho_1}\|^2_2-\|W_{\rho_0}\|^2_2|&\geq &\frac{1}{(2\pi)^2}|\int_{\mathbb{R}^2}|\widetilde{V}_{\tilde{\delta}}(w)|^2dw|
=\frac{1}{(2\pi)^2}|\pi\int_{\mathbb{R}}|t||J_{\tilde{\delta}}(t)|^2dt|\\
&\geq &\frac{1}{(2\pi)^2}2^2\pi^2 L\alpha r\tilde{\delta}^{2-r}e^{2\alpha/\tilde{\delta}^r}2
\int_{(2\lambda+\frac{1}{\tilde{\delta}^r})^{1/r}}^{(D-2\lambda+\frac{1}{\tilde{\delta}^r})^{1/r}}t e^{-4\alpha t^r}dt\\
&=&2 L\alpha r\tilde{\delta}^{2-r}e^{2\alpha/\tilde{\delta}^r} \int_{(2\lambda+\frac{1}{\tilde{\delta}^r})^{1/r}}^{(D-2\lambda+\frac{1}{\tilde{\delta}^r})^{1/r}}t e^{-4\alpha t^r}dt\\
&\geq &\frac{1}{2}Le^{2\alpha /\tilde{\delta}^r}\left(e^{-4\alpha(2\lambda+\frac{1}{\tilde{\delta}})}\left(1+o(1)\right)-e^{-4\alpha
(D-2\lambda+\frac{1}{\tilde{\delta}})}\left(1+o(1)\right)\right)\\
&\geq &\frac{1}{2}Le^{-2\alpha/\tilde{\delta}^r}\left(e^{-8\alpha\lambda}-e^{-4\alpha(D-2\lambda)}\right)\left(1+o(1)\right)\\
&=&2\phi_n\left(e^{-8\alpha\lambda}-e^{-4\alpha(D-2\lambda)}\right)\left(1+o(1)\right)
\end{eqnarray*}
for $n$ large enough, with $\phi_n=\frac{1}{4}\varphi_n$ where $\varphi_n$ is the rate of convergence define in \eqref{vitesse}. Note that we obtain lower bounds for $\tilde{\delta}$ solution of \eqref{equah} for the case $0<r<2$, while we obtain optimal rates (up to a logarithmic factor) of order $(n\log n)^{-\frac{\alpha}{a+\alpha}}$ for $r=2$, with $a$ defined in \eqref{equahdeux}.
\paragraph{\textbf{Proof of \eqref{hyp3}}}
Let us now bound $n\chi^2$. From the lemma 6.1 we get that $p_0(x)\geq Cx^{-2}$ for all $|x|\geq 1$. After a convolution with the gaussian density of the noise the asymptotic decay can not be faster $$p_0^\eta(y)\geq \frac{C_1}{y^2},\forall|y|\geq M,$$ for some fixed $M>0$. Notice that $C$ design a constant which may change along the proof.
\begin{eqnarray}
n\chi^2&\leq&\pi\int\frac{(p_1^\eta(y)-p_0^\eta(y))^2}{p_0^\eta(y)}dy\nonumber\\
\label{equaun}
&\leq &Cn\left(C(M)\|p_1^\eta(y)-p_0^\eta(y)\|^2_2+\int_{|y|>M}y^2\left(p_1^\eta(y)-p_0^\eta(y)\right)^2dy\right).
\end{eqnarray}
In the first term we have
\begin{eqnarray}
\|p_1^\eta(y)-p_0^\eta(y)\|^2_2&=&C\int|J_{\tilde{\delta}}(t)|^2e^{-(1-\eta)t^2/(2\eta)}dt\nonumber\\
&\leq &C\tilde{\delta}^{2-r}e^{2\alpha/\tilde{\delta}^r}\int^{(D-\lambda+\frac{1}{\tilde{\delta}^r})^{1/r}}_{(\lambda+\frac{1}{\tilde{\delta}^r})^{1/r}}
e^{-(1-\eta)t^2/(2\eta)-4\alpha t^r}dt\nonumber\\
&\leq &C\tilde{\delta}^{3-r}e^{2\alpha/\tilde{\delta}^r}\int^{\infty}_{(\lambda+\frac{1}{\tilde{\delta}^r})^{1/r}}
te^{-(1-\eta)t^2/(2\eta)-4\alpha t^r}dt\nonumber\\
\label{equadeux}
&\leq &C\tilde{\delta}^{3-r}\exp{\left(-\frac{2\alpha}{\tilde{\delta}^r}-\frac{1-\eta}{2\eta\tilde{\delta}^2}\right)}.
\end{eqnarray}
Let us see the second part of the sum
\begin{eqnarray}
\int_{|y|>M} y^2(p_1^\eta(y)-p_0^\eta(y))^2dy&\leq &\int(\frac{\partial}{\partial 
t}(J_{\tilde{\delta}}(t)e^{-(1-\eta)t^2/(4\eta)}))^2dt\nonumber\\
&\leq &C\tilde{\delta}^{2-r}e^{2\alpha/\tilde{\delta}^r}\int^{\infty}_{(\lambda+\frac{1}{\tilde{\delta}^r})^{1/r}}t^2e^{-(1-\eta)t^2/(2\eta)}e^{-4\alpha t^r}dt\nonumber\\
\label{equatrois}
&\leq 
&C\tilde{\delta}^{1-r}\exp{\left(-\frac{2\alpha}{\tilde{\delta}^r}-\frac{1-\eta}{2\eta\tilde{\delta}^2}\right)}.
\end{eqnarray}
In case $0<r<2$, by taking $\tilde{\delta}$ as solution of \eqref{equah} we have the expressions in \eqref{equadeux} and \eqref{equatrois} tend to $0$ and together with \eqref{equaun} conclude. For the case $r=2$, we proved a weaker form of \eqref{hyp3}: $n\chi^2=O(1)$. As $\tilde{\delta}$ given in \eqref{equahdeux}, we have the expression in \eqref{equadeux} tend to $0$ while the expression in \eqref{equatrois} stays bounded as $n\rightarrow\infty$ and together with \eqref{equaun} we get the wanted result.
\bibliographystyle{plain}
\bibliography{bibkatia}

\end{document}